\documentclass[12pt]{article}
\usepackage{amssymb,latexsym,amsmath,graphics}
\usepackage{mathrsfs} 

\newtheorem{defn}{Definition}[section]

\def\endproof{\relax\ifmmode\expandafter\endproofmath\else
 \unskip\nobreak\hfil\penalty50\hskip.75em\hbox{}\nobreak\hfil\bull
  {\parfillskip=0pt \finalhyphendemerits=0 \bigbreak}\fi}

\def\bull{\vbox{\hrule\hbox{\vrule
    \kern3pt\vbox{\kern6pt}\kern3pt\vrule}\hrule}}

\def\noi{\noindent}

\def\dsp{\displaystyle}

\def\S{\sum \limits_{n= 0}^\infty}

\def\C{\mathscr{C}}
\def\A{\ifmmode{\alpha_{ij}}\else{${\alpha_{ij}}$}\fi}
\def\B{\ifmmode{\beta_{ij}}\else{${\beta_{ij}}$}\fi}
\def\Al{\ifmmode{\alpha}\else{${\alpha}$}\fi}
\def\Be{\ifmmode{\beta}\else{${\beta}$}\fi}
\def\Si{\ifmmode{\sigma}\else{${\sigma\,}$}\fi}
\def\Ph{\ifmmode{\Phi}\else{${\Phi\,}$}\fi}
\def\X{\ifmmode{x_{ij}}\else{${x_{ij}\,}$}\fi}
\def\N{\ifmmode{n_{ij}}\else{${n_{ij}\,}$}\fi}

\DeclareMathOperator{\cyc}{cyc}
\DeclareMathOperator{\fix}{fix}

\def\Sym{\mathfrak{S}}

\title{Multilinear generating functions for Charlier polynomials}

\author{Ira M. Gessel\thanks
   {Partially supported by NSF Grant  DMS-0200596}\\
\small Department of Mathematics\\[-0.8ex]
\small Brandeis University, Waltham, MA, USA\\[-0.8ex]
\small \texttt{gessel@brandeis.edu}
\and
Pallavi Jayawant\\
\small Department of Mathematics\\[-0.8ex]
\small Bates College, Lewiston, ME, USA\\[-0.8ex]
\small \texttt{pjayawan@bates.edu}}

\date{\small Mathematics Subject Classifications: 05A15, 05A19, 05A40, 33C45}

\begin{document}

\maketitle

\begin{abstract} Charlier configurations provide a combinatorial model for Charlier polynomials. We use this model to give a combinatorial proof of a multilinear generating function for Charlier polynomials. As special cases of the multilinear generating function, we obtain the bilinear generating function for Charlier polynomials and formulas for derangements. 
\end{abstract}

\section{Introduction}\label{S:intro}

Charlier polynomials have been studied using combinatorial methods in \cite{foata}, \cite{gessel1}, \cite{labelle-yeh}, \cite{labelle-yeh1}, \cite{viennot}, and  \cite{zeng}. In this paper, we prove a multilinear generating function for Charlier polynomials using the combinatorial model of Charlier configurations \cite{labelle-yeh,labelle-yeh1} and the approach of Foata and Garsia \cite{foata-garsia:multi} in their proof of Slepian's multilinear extension of the Mehler formula for Hermite polynomials \cite{slepian}. We then obtain some formulas for derangements as special cases of this generating function.

The Charlier polynomials are usually defined by the formula
\[c_n(a,r) = \,_2F_0(-n,-a;-;-r^{-1})=\sum_{k=0}^n \binom nk (-a)_k \frac{r^{-k}}{k!},\]
where $(u)_k = u(u+1)\cdots(u+k-1)$.
In order to assign convenient weights in the combinatorial model, we work with renormalized Charlier polynomials $C_n(a,r)$ defined by 
$$C_n(a,r) = r^nc_n(-a,r) =\sum\limits_{k=0}^n {n\choose k} (a)_k r^{n-k} .$$ 
Our main result is the multilinear generating function
\begin{multline}\label{M}\sum_{(n_{ij})} \frac
{\dsp{\prod_{1 \le i < j \le k} {\X^{n_{ij}} }}}{
\dsp{\prod_{1 \le i < j \le k}n_{ij}!}} C_{n_1}(a_1,r_1) \cdots
C_{n_k}(a_k,r_k) \\[-30pt]
=\prod_{1 \le i < j \le k} e^{r_i r_j \X }\,
\sum_{(n_{ij})} \prod_{1 \le i \le k}\dsp{\frac{(a_i)_{n_i}}{{(1-\sum_{j\ne i} r_j\X)}^{n_i+a_i}}}
\frac
{\dsp{\prod_{1 \le i < j \le k} {{x_{ij}}^{n_{ij}} }}}{
\dsp{\prod_{1 \le i < j \le k}n_{ij}!}} \, , 
\end{multline}
where each sum runs over all $k\times k$ symmetric matrices $(\N)$ with non-negative integral
entries and with diagonal entries zero, $n_i=\sum\limits_{j=1}^k \N$ for $1\le i\le k$, and $x_{ij}=x_{ji}$.

To give a combinatorial proof of (\ref{M}), we begin with a discussion of Charlier configurations and their representation by digraphs in section \ref{S:Ch}. Then in section \ref{S:Ml} we give the combinatorial proof of the multilinear generating function. The main idea of the proof is to show that both sides of the formula count the same set of digraphs. We discuss the special cases of the multilinear generating function in section \ref{S:relformulas}.   


\section{Charlier Configurations} \label{S:Ch}

Let $[n]$ denote the set $\{1,2,\ldots, n\}$.

\begin{defn}
 A Charlier configuration on the set $S$ is a pair 
$\Phi = ((A, \sigma), B)$, where $(A,B)$ is an ordered partition of $S$ and $\sigma$ is a
permutation of $A$. 
\end{defn}

\begin{figure}[h]
\begin{center}
\scalebox{0.8}{\includegraphics{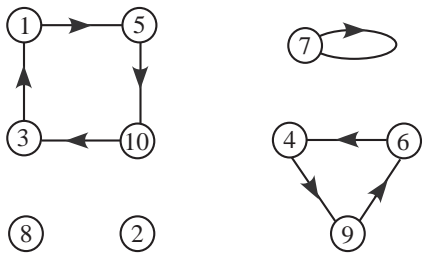}}
\caption{Charlier configuration on $[10]$}\label{Fi:charconfig}
\end{center}
\end{figure}

A Charlier configuration is called a partial permutation in \cite{labelle-yeh}. The configuration $\Ph$ can be represented by a digraph with
vertex set $S$ and with an edge from $i$ to $j$ if and only if $\sigma (i)=j$. 
Figure \ref{Fi:charconfig} shows a Charlier configuration on $[10]$.

Here $A=\{1, 3, 4, 5, 6, 7, 9, 10\}$, $B=\{2, 8\}$, and 
$\sigma=(7)\,\, (4\,\, 9\,\, 6)\,\, (1\,\, 5\,\, 10\,\, 3)$  in disjoint cycle notation.

\subsection{Combinatorial interpretation of Charlier polynomials.}\label{combint}
We assign a
weight to a Charlier configuration $\Ph$ by  assigning a weight $a$ to each cycle of $\sigma$ and a weight $r$ to
each point of $B$. Then the weight of $\Ph$ is
$a^{\cyc(\sigma)}r^{|B|}$ where $\cyc(\sigma)$ denotes the number of cycles in $\sigma$. (If we had not renormalized the Charlier polynomials, we would assign the weight $-a$ to each cycle of $\sigma$ and the weight $1/r$ to each point of $A$.) 
We use the following well-known facts about generating functions for permutations, which are proved, for example,  in \cite[p.~19]{stanley:vol1}.

\noi {\bf Fact 1.} $\sum_{k=0}^n c(n,k)a^k=(a)_n$ where $c(n,k)$ is the number of permutations of $[n]$ with exactly $k$ cycles (the unsigned Stirling number of the first kind).

\noi {\bf Fact 2.} The exponential generating function for all permutations, with cycles weighted by $a$, is $(1-z)^{-a}$.

Let $\C_S$ denote the set of Charlier configurations on $S$. Then it follows easily from Fact~1 that   $C_n(a, r)$ is sum of the weights of the elements of $\C_{[n]}$. 
 
\section{Combinatorial Proof of the Multilinear Formula} \label{S:Ml}

We assume that the reader is familiar with enumerative
applications of exponential generating functions, as described,
for example, in \cite[Chapter 5]{stanley:vol2} and \cite{bergeron}. 
The product formula and the exponential formula for exponential
generating functions discussed in these references play an
important role in the combinatorial proof of the multilinear formula. The theory of species (as used in \cite{labelle-yeh} and \cite{labelle-yeh1}) could be used to provide a proof of the formula as well. 

The formula \eqref{M} could be proved by interpreting it as a multivariable exponential generating function in the variables $x_{ij}$, which would require the use of digraphs with multiple sets of labels. The proof is simpler if we use exponential generating functions in only one variable, so that we can use a single set of labels.
To accomplish this, we rewrite the formula by replacing \X  with $z\X $. Now we can think of the
formula as an exponential generating function in the single variable $z$. The formula  is now
\begin{multline}\label{mlinear}\sum_{(n_{ij})} \frac
{\dsp{\prod_{1 \le i < j \le k} {\X^{n_{ij}} }}}{
\dsp{\prod_{1 \le i < j \le k}n_{ij}!}} C_{n_1}(a_1,r_1) \cdots
C_{n_k}(a_k,r_k) z^{\sum \N}=\\ 
\prod_{1 \le i < j \le k}\! \! \! \! e^{r_i r_j \X z}\! \!
\prod_{1 \le i \le k}\dsp{\frac{1}{{(1-z\sum_{j\ne i} r_j\X)}^{a_i}}}
\sum_{(n_{ij})} \prod_{1 \le i \le k}\dsp{\frac{(a_i)_{n_i}}{{(1-z\sum_{j\ne i} r_j\X)}^{n_i}}}
\frac
{\dsp{\prod_{1 \le i < j \le k} {\X^{n_{ij}} }}}{
\dsp{\prod_{1 \le i < j \le k}n_{ij}!}} z^{\sum \N} \, . \end{multline}
We will prove this formula, which is equivalent to \eqref{M}. We begin with a description of the digraphs counted by the left side of the formula.

\subsection{Digraphs counted by the left side.}
We can rewrite the left side of \eqref{mlinear} as follows:
$$\sum_{n\ge 0} \frac{z^n}{n!}
\sum_{\substack{(\N) \\ \sum\limits_{i<j}\N=n}}\frac{n!} {\prod\limits_{1\le i<j\le k}
\N!} 
\prod_{1 \le i < j \le k} {\X^{n_{ij}} }
C_{n_1}(a_1,r_1) \cdots
C_{n_k}(a_k,r_k) \, .$$

Let $(N_{ij})_{1\le i<j\le k}$ be an ordered partition of $[n]$ such that
$|N_{ij}|=\N$. 
For $j>i$, let $N_{ji}=N_{ij}$. Let $N_i=\cup_{j\ne i}N_{ij}$. Then $N_i\cap N_j=N_{ij}$. Since $n_i=\sum_j \N$,  it follows that $|N_i|=n_i$.
Let $H$ be the set of all ordered tuples $((N_{ij}), \Ph_1, \ldots ,\Ph_k)$ such that
\begin{enumerate}
\item $(N_{ij})$ is an ordered partition of $[n]$ with the above properties.

\item Each $\Ph_i$ is a Charlier configuration on $N_i$, i.e., $\Ph_i \in \C_{N_i}$.
\end{enumerate}
Then each point of $[n]$ is in exactly two configurations. This follows from the
fact that each point is in exactly one $N_{ij}$ and $N_i\cap N_j=N_{ij}$.
To the Charlier configuration $\Phi_i=((A_i, \sigma_i), B_i)$ we assign the weight
$a_i^{\cyc(\sigma)}r_i^{|B|}$. We also assign an additional weight of $\X$ to each point of $N_{ij}$. The weight of a tuple in $H$ is defined to be the product of the weights of its constituent Charlier configurations and its points. Then it is easy to see that the left side of 
\eqref{mlinear} is the exponential 
generating function for $H$ with these weights. 

\begin{figure}[h]
\begin{center}
\scalebox{0.8}{\includegraphics{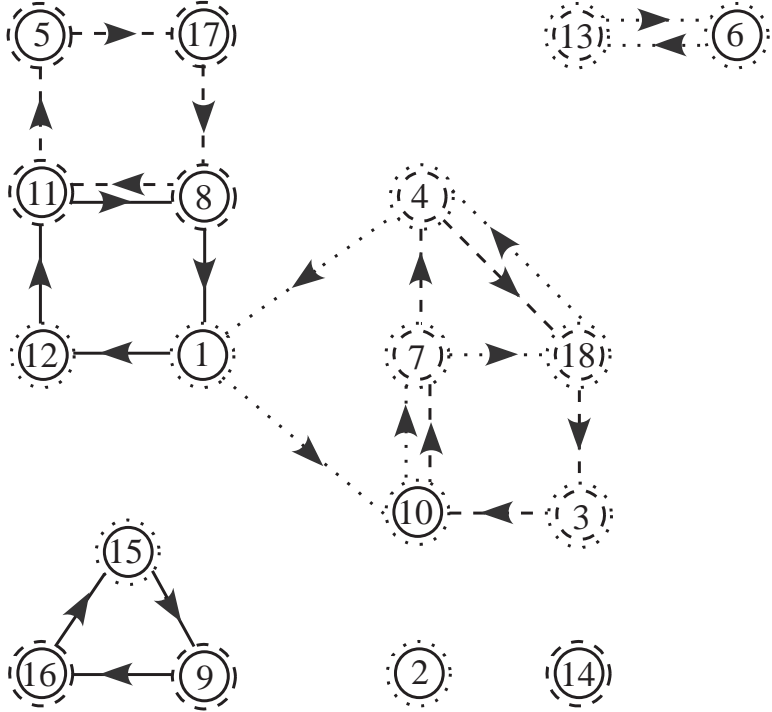}}
\caption{Digraph counted by left side for $k=3$}\label{Fi:trigraph}
\end{center}
\end{figure}

We associate a digraph to a tuple $((N_{ij}), \Ph_1, \ldots ,\Ph_k)$  in $H$
by  superimposing the digraphs of these $k$ Charlier configurations on $[n]$ in which each $\Ph_i$ is on $n_i$ of these vertices. Figure \ref{Fi:trigraph} shows such a digraph for $k=3$.
The configurations $\Ph_1, \Ph_2, \Ph_3$ are respectively represented by solid lines, dashed lines, and dotted lines.
Each vertex is in exactly two configurations and this is indicated by the two different circles around each vertex. 

The tuple and the configurations corresponding to Figure \ref{Fi:trigraph} are given by
\begin{align*}
N_{12}&=\{5, 8, 9, 11, 14, 16, 17\},
 N_{13}=\{1, 2, 6, 12, 15\}, \text{ and }
 N_{23}=\{3, 4, 7, 10, 13, 18\}.\\
\Phi_1&=((A_1, \sigma_1), B_1), \text{ where}\\
     &\quad A_1=\{1, 8, 9, 11, 12, 15, 16\},
      B_1=\{2, 5, 6, 14, 17\},  \\
     &\quad  \sigma_1=(9\,\, 16\,\, 15)\,\, (1\,\, 12\,\, 11\,\, 8).\\
\Phi_2&=((A_2, \sigma_2), B_2), \text{ where}\\
    &\quad A_2=\{3, 4, 5, 7, 8, 10, 11, 17, 18\},
     B_2=\{9, 13, 14, 16\},\\
    &\quad\sigma_2=(5\,\, 17\,\, 8\,\, 11)\,\, (3\,\, 10\,\, 7\,\, 4\,\, 18).\\
\Phi_3&=((A_3, \sigma_3), B_3), \text{ where}\\
  &\quad A_3=\{1, 4, 6, 7, 10, 13, 18\}, B_3=\{2, 3, 12, 15\},\\
     &\quad \sigma_3=(6\,\, 13)\,\, (1\,\, 10\,\, 7\,\, 18\,\, 4).
\end{align*}

We may identify $H$
with the set of these digraphs, for which  so that the left side  of \eqref{mlinear} is the exponential generating function.
We now enumerate these digraphs in another way: We consider  their connected components, which are of three types, and use the product formula for exponential generating functions to show that the right side of \eqref{mlinear} is also a generating function for $H$.

\subsection{Connected components of digraphs in $H$.}
Let $((N_{ij}), \Ph_1, \ldots ,\Ph_k)$ be a tuple in $H$, where  
$\Ph_i=((A_i, \sigma_i), B_i)\in C_{n_i}$, for $1\le i\le k$. The connected components of the digraph representing this
tuple are  of the following three types:

For $i<j$, a {\bf type $1_{ij}$ connected component} is an isolated vertex which is in $\Ph_i$ and $\Ph_j$ but not in $\sigma_i$ or $\sigma_j$.
In Figure \ref{Fi:trigraph}, vertex $2$ is of type $1_{13}$ and vertex $14$  is  of type $1_{12}$.
Such a vertex belongs to $B_i\cap B_j$ and is weighted  by $r_i r_j \X$.
It follows that the  exponential generating function for digraphs all of whose components are of type 
$1_{ij}$, which  we call {\bf type $1_{ij}$ digraphs}, is $e^{r_i r_j \X z}$.

A {\bf type $2_i$ connected component} is a cycle of $\sigma_i$ in which no vertex is in any other~$\sigma_j$. In Figure \ref{Fi:trigraph},  the cycle $(6\,\, 13)$ is a type $2_3$ component and $(9\,\, 16\,\, 15)$ is a type~$2_1$ component.
The cycle of a type $2_i$ component 
weighted by $a_i$ and each vertex of the cycle is in some $B_j$ and so is weighted by $r_j$ and \X. 
A {\bf type $2_i$ digraph} is a
digraph in which every connected component is of type $2_i$. Such digraphs can be considered as 
permutations in which each cycle is weighted by $a_i$ and each vertex is weighted
by some $r_j$ and $\X$ for some $j$. By a slight modification of Fact 2 in subsection \ref{combint}, it follows that the
exponential generating function for type $2_i$ digraphs is $(1-z\sum\limits_{j\ne i} r_j\X)^{-a_i}$.

Any connected component that is not of type $1_i$ or $2_{ij}$ is called a 
{\bf type 3 connected component}. In a type 3 connected component, every vertex is in at least one
permutation and every cycle contains at least one vertex that is also in another permutation.

\begin{figure}[h]
\begin{center}
\scalebox{0.8}{\includegraphics{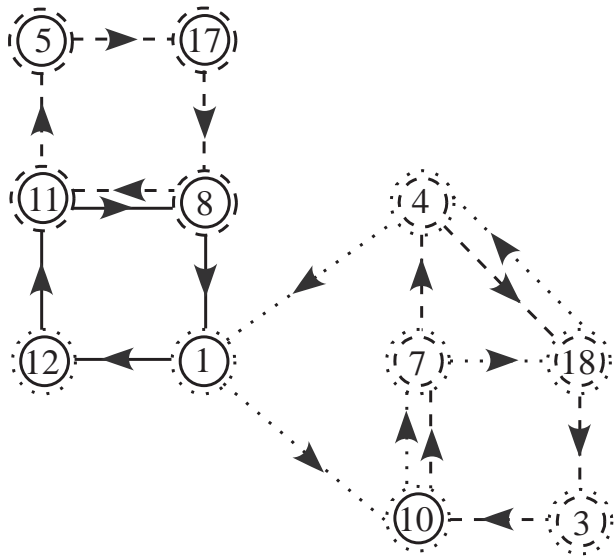}}
\caption{Type 3 connected component}
\label{Fi:bigcomponent}
\end{center}
\end{figure}

\begin{figure}[h]
\begin{center}
\scalebox{0.8}{\includegraphics{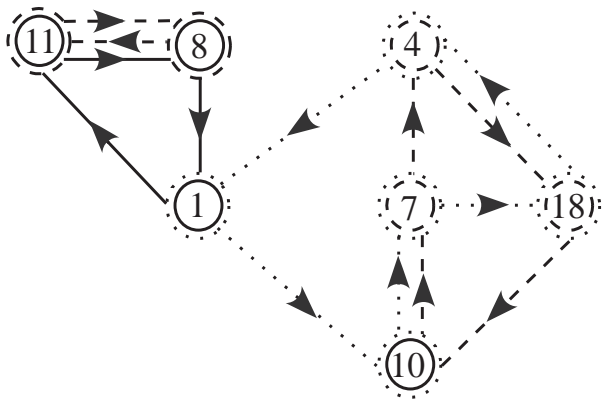}}
\caption{Reduced digraph on 7 vertices}
\label{Fi:reduced}
\end{center}
\end{figure}

The type 3 connected component from Figure \ref{Fi:trigraph}
is shown in Figure \ref{Fi:bigcomponent}.
A {\bf type 3 digraph} is a digraph all of whose connected components are of type 3. We say
that a type 3 digraph is reduced if every vertex is in two permutations. 
Thus, a reduced type 3 
digraph on $n$ vertices is an ordered tuple $((N_{ij}), \Ph_1, \ldots ,\Ph_k)$ 
where each $\Ph_i=((A_i, \sigma_i), \varnothing)$; i.e., each $\Ph_i$ is simply a permutation $\sigma_i$ on 
$n_i$ vertices. Figure \ref{Fi:reduced} shows a reduced type 3 digraph on 7 vertices. 
Since each cycle of $\sigma_i$ is weighted by $a_i$, the exponential generating
function for reduced type 3 digraphs is 
$$\sum_{n\ge 0}\frac{z^n}{n!}\sum_{\substack{(\N) \\ \sum\limits_{i<j}\N=n}}
\frac{n!}{\prod\limits_{1\le i<j\le k}
\N!} \prod_{1\le i\le k}(a_i)_{n_i} \prod_{1\le i<j\le k} \X^{\N}\, .$$
Each vertex in a reduced type 3 digraph has two outgoing edges belonging to two different
permutations. Any type 3 digraph can be obtained
from a reduced type 3 digraph by replacing each outgoing edge at every vertex by a sequence of
ordered edges. An outgoing edge in $\sigma_i$ is replaced by a sequence of edges, such that
each new vertex is in $\sigma_i$, but not in any other $\sigma_j$. 
Hence each new vertex is weighted by
$r_j$ and $\X$ for some $j$.  Thus the exponential generating function for type 3 digraphs is
$$\sum_{n\ge 0}{z^n\over n!}\sum_{(\N)\atop \sum\limits_{i<j}\N=n}{n!\over \prod\limits_{1\le i<j\le k}
\N!} \prod_{1\le i\le k}{(a_i)_{n_i}\over (1-z\sum\limits_{j\ne i} r_j\X)^{n_i}}
 \prod_{1\le i<j\le k} \X^{\N}\, .$$

It follows from the product formula for exponential generating functions that the exponential generating function for 
digraphs in $H$ is the product of the generating functions for all of the types of digraphs described above and this is
\begin{multline*}\prod_{1 \le i < j \le k} e^{r_i r_j \X z}
\prod_{1\le i\le k}\frac 1 {(1-z\sum\limits_{j\ne i} r_j\X)^{a_i}}\\ 
\times\sum_{n\ge 0}{z^n\over n!}\sum_{(\N)\atop \sum\limits_{i<j}\N=n}{n!\over \prod\limits_{1\le i<j\le k}
\N!} \prod_{1\le i\le k}{(a_i)_{n_i}\over (1-z\sum\limits_{j\ne i} r_j\X)^{n_i}}
 \prod_{1\le i<j\le k} \X^{\N}\, ,
 \end{multline*}
which is equal to the right side of \eqref{mlinear}.


\section{Specializations}\label{S:relformulas}

For $k=2$, the only parameter $n_{ij}$ in \eqref{M} is $n_{12}$, and $n_1=n_2=n_{12}$.
If we write $n$ for $n_{12}$, $a$ for $a_1$, $b$ for $a_2$, $r$ for $r_1$, $s$ for $r_2$, and $x$ for $x_{12}$ then 
the multilinear formula \eqref{M} reduces to the bilinear formula
\begin{equation}\label{B} 
\sum_{n\ge 0} C_n(a, r) C_n(b, s){x^n\over n!}=e^{rsx}
\sum_{n\ge 0} {(a)_n (b)_n\over (1-sx)^{n+a}\, (1-rx)^{n+b}}{x^n\over n!}\, . \end{equation}
Formula (2.47) in Askey's book \cite{askey} is equivalent to the case of \eqref{B} in which $a$ and $b$ are negative integers, and the general case is easily derived from this. Note that $a$ and $b$ are switched on the right side of the formula in the book.

Similarly, the case $k=3$ of \eqref{M} may be written
\begin{multline}
\label{k=3}
\sum_{l,m,n} C_{l+m}(a,r) C_{l+n}(b,s) C_{m+n}(c,t) \frac{x^l}{l!}\frac{y^m}{m!}\frac{z^n}{n!}= e^{rsx+rty+stz}\\
 \times\sum_{l,m,n}
 \frac{(a)_{l+m} (b)_{l+n} (c)_{m+n}}{(1-sx-ty)^{l+m+a} (1-rx-tz)^{l+n+b} (1-ry-sz)^{m+n+c}}
 \frac{x^l}{l!}\frac{y^m}{m!}\frac{z^n}{n!}.
\end{multline}

Some special  cases of these formulas are worth mentioning. Setting $x=0$ in \eqref{k=3} gives
\begin{multline}
\label{mspecial}
\qquad\sum_{m,n} C_{m}(a,r) C_{n}(b,s) C_{m+n}(c,t) \frac{y^m}{m!}\frac{z^n}{n!}\\
 = e^{(ry+sz)t}\sum_{m,n}
 \frac{(a)_{m} (b)_{n} (c)_{m+n}}{(1-ty)^{m+a} (1-tz)^{n+b} (1-ry-sz)^{m+n+c}}
   \frac{y^m}{m!}\frac{z^n}{n!}.\qquad
\end{multline}
Formula \eqref{mspecial} may be viewed as a Charlier polynomial analogue of a formula of Carlitz~\cite{carlitz:mehlerext} for Hermite polynomials, which is  a special cases of Slepian's
multilinear extension of the Mehler formula \cite{slepian}.

By applying the fact that $C_n(0,1)=1$, we can find other simplifications. Thus setting $b=0$ and $s=1$ in \eqref{B} gives the usual exponential generating function for Charlier polynomials,
\begin{equation*}
\sum_{n\ge 0} C_n(a, r) {x^n\over n!}=e^{rx} (1-x)^{-a}\,.
\end{equation*}
Setting $c=0$ and $t=1$ in \eqref{k=3} gives a generalization of \eqref{B}:
\begin{equation}
\label{M3}
\sum_{l,m,n} C_{l+m}(a,r) C_{l+n}(b,s) \frac{x^l}{l!}\frac{y^m}{m!}\frac{z^n}{n!}
 =e^{rsx+ry+sz}\sum_{l}
 \frac{(a)_{l} (b)_{l}}{(1-sx-y)^{l+a} (1-rx-z)^{l+b}}
 \frac{x^l}{l!}.
\end{equation}

\subsection{Permutations.}
\label{S:specialbilinear}

The Charlier polynomials can be normalized in another way so as to count permutations by cycles and fixed points.
We define polynomials $D_n(\alpha, u)$ by 
\[D_n(\alpha, u) = \sum_{\pi\in \Sym_n}\alpha^{\cyc_{>1} (\pi)}u^{\fix(\pi)},
\]
where $\Sym_n$ is the set of permutations of $[n]$,
$\cyc_{>1}(\pi)$ is the number of cycles of $\pi$ of length greater than 1, 
and $\fix(\pi)$ is the number of fixed points of $\pi$. 

We can express the polynomials $D_n(\alpha, u)$ in terms of Charlier polynomials.
To a Charlier configuration $\Phi = ((A, \sigma), B)$ on $[n]$ we may associate the permutation $\pi$ of $[n]$
such that $\pi(i) =\sigma(i)$  for $i\in A$ and $\pi(i)=i$ for $i\in B$. Conversely, given a permutation $\pi$ of 
$[n]$, the corresponding Charlier configurations may be constructed by choosing an arbitrary subset $B$ of the set of fixed points of $\pi$ and taking $\sigma$ to be the restriction of $\pi$ to $A=[n]\setminus B$.  
This construction yields the relation
\[C_n(a, r) = D_n(a, a+r)\]
and thus
\[D_n(\alpha, u) = C_n(\alpha, u-\alpha).\]
So formulas \eqref{B}--\eqref{M3} may rewritten as generating functions for the polynomials $D_n(\alpha, u)$.
Of particular interest are the specializations $D_n(\alpha) = D_n(\alpha,0) = C_n(\alpha, -\alpha)$, which count derangements (permutations without fixed points) by cycles, and $D_n=D_n(1)=C_n(1,-1)$, the number of derangements of $[n]$.

Making the appropriate substitutions in \eqref{B} gives 
\begin{equation}\label{derange}
\S D_n(\alpha) D_n(\Be) \frac {x^n}{n!} = e^{\alpha\Be x} \S \frac {(\alpha)_n\, (\Be)_n}{(1+\Be x)^{n+\alpha}\, (1+\alpha
x)^{n+\Be}}\frac{x^n}{n!}.
\end{equation}
Formula \eqref{derange} was proved by Gessel \cite{gessel} as a special case of a  generating function
for $3 \times n$ Latin rectangles, using an approach similar to that of this paper. (See also \cite{agj} and \cite{zeng1}.)

The corresponding specialization of \eqref{mspecial} is 
\begin{multline*}
\sum_{m,n=0}^{\infty} D_m(\alpha) D_n(\beta) D_{m+n}(\gamma)\frac {y^m}{m!}\frac{z^n}{n!} \\ 
=e^{(\alpha y+\beta z)\gamma} \sum_{m,n=0}^{\infty} \frac {(\alpha)_m(\beta)_n(\gamma)_{m+n}}
{(1+\gamma y)^{m+\alpha} (1+\gamma z)^{n+\beta} (1+\alpha y+\beta z)^{m+n+\gamma}}
\frac {y^m}{m!}\frac{z^n}{n!}
\end{multline*}
and of \eqref{M3} is
\begin{equation*}
\sum_{l,m,n} D_{l+m}(\alpha) D_{l+n}(\beta) \frac{x^l}{l!}\frac{y^m}{m!}\frac{z^n}{n!}
 =e^{\alpha\beta x-\alpha y-\beta z}\sum_{l}
 \frac{(\alpha)_{l} (\beta)_{l}}{(1+\beta x-y)^{l+\alpha} (1+\alpha x-z)^{l+\beta}}
 \frac{x^l}{l!}.
\end{equation*}



\begin{thebibliography}{9}

\bibitem{agj}
G. E. Andrews, I. P. Goulden, and D. M. Jackson, Generalizations of Cauchy's summation theorem for Schur functions,
Trans. Amer. Math. Soc. \textbf{310} (1988), 805--820.

\bibitem{askey}
R. Askey, \emph{Orthogonal polynomials and special functions}, Society for
  Industrial and Applied Mathematics, Philadelphia, Pa., 1975.

\bibitem{bergeron}
F. Bergeron, G. Labelle, and P. Leroux, \emph{Combinatorial
Species and Tree-like Structures}, 
  Encyclopedia of Mathematics and
  its Applications, Vol.~67, Cambridge University Press, Cambridge, 1997.
  Translated
  from the 1994 French original by Margaret Readdy.

\bibitem{carlitz:mehlerext}
L.~Carlitz, \emph{Some extensions of the {M}ehler formula}, Collect. Math.
  \textbf{21} (1970), 117--130.

\bibitem{foata} D. Foata, \emph{Combinatoire des identit\'es sur les polyn\^omes orthogonaux}, 
Proceedings of the International Congress of Mathematicians (Warsaw, 1983), ed.
 Zbigniew Ciesielski and Czes\l aw Olech,
PWN---Polish Scientific Publishers, Warsaw; North-Holland Publishing Co., Amsterdam,
1984, pp. 1541--1553.

\bibitem{foata-garsia:multi}
D. Foata and A.~M. Garsia, \emph{A combinatorial approach to the
  {M}ehler formulas for {H}ermite polynomials},  Relations between combinatorics
  and other parts of mathematics (Proc. Sympos. Pure Math., Ohio State Univ.,
  Columbus, Ohio, 1978), ed. D. K. Ray-Chaudhuri, Amer. Math. Soc., Providence, R.I., 1979,
 pp. 163--179.

\bibitem{gessel}
I. M. Gessel, \emph{Counting three-line {L}atin rectangles}, Combinatoire
  \'enum\'erative (Montreal, Que., 1985/Quebec, Que., 1985), ed. \ G. Labelle and P. Leroux,  Lecture notes in
  Math.~1234, Springer, Berlin, 1986, pp. 106--111.

\bibitem {gessel1}
I. M. Gessel, \emph{Generalized rook polynomials and orthogonal polynomials}, $q$-Series and Partitions, ed. D.  Stanton,  IMA Volumes in Math. and its Appl. \textbf{18}, Springer-Verlag, New York, 1989, pp. 159--176.

\bibitem{kibble}
W.~F. Kibble, \emph{An extension of a theorem of {M}ehler's on {H}ermite
  polynomials}, Proc. Cambridge Philos. Soc. \textbf{41} (1945), 12--15.

\bibitem{labelle-yeh} J. Labelle and Y. N. Yeh, \emph{The combinatorics of Laguerre, Charlier, and Hermite polynomials}, Stud. in Appl. Math. \textbf{80} (1989), 25--36.

\bibitem{labelle-yeh1} J. Labelle and Y. N. Yeh, \emph{Combinatorial proofs of some limit formulas involving orthogonal polynomials}, Discrete Math. \textbf{79} (1989/90), 77--93.

\bibitem{slepian}
D. Slepian, \emph{On the symmetrized {K}ronecker power of a matrix and
  extensions of {M}ehler's formula for {H}ermite polynomials}, SIAM J. Math.
  Anal. \textbf{3} (1972), 606--616.

\bibitem{stanley:vol1}
R. P. Stanley, \emph{Enumerative combinatorics, {V}ol. 1}, Cambridge
  University Press, Cambridge, 1997.

\bibitem{stanley:vol2}
R. P. Stanley, \emph{Enumerative combinatorics, {V}ol. 2}, Cambridge University
  Press, Cambridge, 1999.

\bibitem{viennot}
X. G. Viennot, \emph{Une th\'eorie combinatoire des polyn\^omes orthogonaux}, Lecture Notes, Publications du LACIM, UQAM, Montr\'eal, 1983. 


\bibitem{zeng}
J. Zeng, \emph{Lin\'earisation de produits de polyn™mes de Meixner, Krawtchouk, et Charlier}, SIAM J. Math. Anal. \textbf{21} (1990), 1349--1368.

\bibitem{zeng1}
J. Zeng, \emph{Counting a pair of permutations and the linearization coefficients for Jacobi
polynomials}, 
Atelier de combinatoire franco-qu\'ebecois, ed. J. Labelle and J.-G. Penaud, Publications 
du LACIM, vol. 10, Universit\'e du Qu\'ebec \`a Montr\'eal, 1992, pp.~243--257.



\end{thebibliography}
\end{document}